\documentclass[reqno]{amsart}
\usepackage{setspace,amssymb}
\usepackage{ifpdf}
\ifpdf
 \usepackage[hyperindex]{hyperref}
\else
 \expandafter\ifx\csname dvipdfm\endcsname\relax
 \usepackage[hypertex,hyperindex]{hyperref}
 \else
 \usepackage[dvipdfm,hyperindex]{hyperref}
 \fi
\fi
\allowdisplaybreaks[4]
\numberwithin{equation}{section}
\theoremstyle{plain}
\newtheorem{theorem}{Theorem}[section]
\newtheorem{prop}{Proposition}[section]
\theoremstyle{remark}
\newtheorem{rem}{Remark}[section]
\theoremstyle{definition}

\DeclareMathOperator{\td}{d\mspace{-1mu}}

\begin{document}

\title[Complete monotonicity and modified Bessel function]
{Complete monotonicity of a difference between the exponential and trigamma functions}

\author[F. Qi]{Feng Qi}
\address[Feng Qi]{Department of Mathematics, School of Science, Tianjin Polytechnic University, Tianjin City, 300387, China; School of Mathematics and Informatics, Henan Polytechnic University, Jiaozuo City, Henan Province, 454010, China}
\email{\href{mailto: F. Qi <qifeng618@gmail.com>}{qifeng618@gmail.com}, \href{mailto: F. Qi <qifeng618@hotmail.com>}{qifeng618@hotmail.com}, \href{mailto: F. Qi <qifeng618@qq.com>}{qifeng618@qq.com}} \urladdr{\url{http://qifeng618.wordpress.com}}

\author[X.-J. Zhang]{Xiao-Jing Zhang}
\address[Xiao-Jing Zhang]{Department of Mathematics, School of Science, Tianjin Polytechnic University, Tianjin City, 300387, China}
\email{\href{mailto: X.-J. Zhang <xiao.jing.zhang@qq.com>}{xiao.jing.zhang@qq.com}}

\begin{abstract}
In the paper, by directly verifying an inequality which gives a lower bound for the first order modified Bessel function of the first kind, the authors supply a new proof for the complete monotonicity of a difference between the exponential function $e^{1/t}$ and the trigamma function $\psi'(t)$ on $(0,\infty)$.
\end{abstract}



\maketitle

\section{Introduction}

In~\cite[Lemma~2]{Yang-Fan-2008-Dec-simp.tex}, the inequality
\begin{equation}\label{e-1-t-1}
\psi'(t)<e^{1/t}-1
\end{equation}
on $(0,\infty)$ was discovered and employed, where $\psi(t)$ denotes the digamma function
\begin{equation*}
\psi(t)=[\ln\Gamma(t)]'=\frac{\Gamma'(t)}{\Gamma(t)}
\end{equation*}
and $\Gamma$ is the classical Euler gamma function which may be defined for $\Re(z)>0$ by
\begin{equation*}
\Gamma(z)=\int^\infty_0t^{z-1} e^{-t}\td t.
\end{equation*}
The functions $\psi'(z)$ and $\psi''(z)$ are respectively called the trigamma function and the tetragamma function. As a whole, the derivatives $\psi^{(k)}(z)$ for $k\in\{0\}\cup\mathbb{N}$ are called polygamma functions.
\par
An infinitely differentiable function $f$ defined on an interval $I$ is said to be a completely monotonic function on $I$ if it satisfies
\begin{equation}
(-1)^kf^{(k)}(x)\ge0
\end{equation}
for all $k\in\{0\}\cup\mathbb{N}$ on $I$. Some properties of the completely monotonic functions can be found in, for example, \cite{psi-proper-fraction-degree-two.tex, widder}.
\par
In~\cite[Theorem~3.1]{QiBerg.tex} and~\cite[Theorem~1.1]{simp-exp-degree-ext.tex}, the following theorem was proved totally by three methods.

\begin{theorem}\label{CM-Exp-thm}
The function
\begin{equation}\label{alpha-exp=psi-eq}
h(t)=e^{1/t}-\psi'(t)
\end{equation}
is completely monotonic on $(0,\infty)$ and
\begin{equation}\label{h(t)-limit=1}
\lim_{t\to\infty}h(t)=1.
\end{equation}
\end{theorem}

The second main result of the paper~\cite{simp-exp-degree-ext.tex} is its Theorem~1.2 which has been referenced in~\cite[Section~1.2]{Bessel-ineq-Dgree-CM.tex} and~\cite[Lemma~2.1]{QiBerg.tex} as follows.

\begin{theorem}\label{exp=k=degree=k+1-int-thm}
For $k\in\{0\}\cup\mathbb{N}$ and $z\ne0$, let
\begin{equation}\label{exp=k=sum-eq-degree=k+1}
H_k(z)=e^{1/z}-\sum_{m=0}^k\frac{1}{m!}\frac1{z^m}.
\end{equation}
For $\Re(z)>0$, the function $H_k(z)$ has the integral representations
\begin{equation}\label{exp=k=degree=k+1-int}
H_k(z)=\frac1{k!(k+1)!}\int_0^\infty {}_1F_2(1;k+1,k+2;t)t^k e^{-zt}\td t
\end{equation}
and
\begin{equation}\label{exp=k=degree=k+1-int-bes}
H_k(z)=\frac1{z^{k+1}}\biggl[\frac1{(k+1)!}+\int_0^\infty \frac{I_{k+2} \bigl(2\sqrt{t}\,\bigr)}{t^{(k+2)/2}} e^{-zt}\td t\biggr],
\end{equation}
where the hypergeometric series
\begin{equation}
{}_pF_q(a_1,\dotsc,a_p;b_1,\dotsc,b_q;x)=\sum_{n=0}^\infty\frac{(a_1)_n\dotsm(a_p)_n} {(b_1)_n\dotsm(b_q)_n}\frac{x^n}{n!}
\end{equation}
for $b_i\notin\{0,-1,-2,\dotsc\}$, the shifted factorial $(a)_0=1$ and
\begin{equation}
(a)_n=a(a+1)\dotsm(a+n-1)
\end{equation}
for $n>0$ and any real or complex number $a$, and the modified Bessel function of the first kind
\begin{equation}\label{I=nu(z)-eq}
I_\nu(z)= \sum_{k=0}^\infty\frac1{k!\Gamma(\nu+k+1)}\biggl(\frac{z}2\biggr)^{2k+\nu}
\end{equation}
for $\nu\in\mathbb{R}$ and $z\in\mathbb{C}$.
\end{theorem}

When $k=0$, the integral representations~\eqref{exp=k=degree=k+1-int} and~\eqref{exp=k=degree=k+1-int-bes} may be written as
\begin{equation}\label{open-answer-1}
e^{1/z}=1+\int_0^\infty \frac{I_1\bigl(2\sqrt{t}\,\bigr)}{\sqrt{t}\,} e^{-zt}\td t
\end{equation}
and
\begin{equation}\label{open-answer-2}
e^{1/z}=1+\frac1{z}\biggl[1+\int_0^\infty \frac{I_{2} \bigl(2\sqrt{t}\,\bigr)}{t} e^{-zt}\td t\biggr]
\end{equation}
for $\Re(z)>0$. Hence, by the well known formula
\begin{equation}\label{polygamma}
\psi^{(n)}(z)=(-1)^{n+1}\int_0^{\infty}\frac{u^n}{1-e^{-u}}e^{-zu}\td u
\end{equation}
for $\Re(z)>0$ and $n\in\mathbb{N}$, see~\cite[p.~260, 6.4.1]{abram}, the function $h(t)$ defined by~\eqref{alpha-exp=psi-eq} has the following integral representation
\begin{equation}\label{h(t)-int-rep-eq}
h(t)=1+\int_0^\infty \biggl[\frac{I_1\bigl(2\sqrt{u}\,\bigr)}{\sqrt{u}\,} -\frac{u}{1-e^{-u}}\biggr]e^{-tu}\td u.
\end{equation}

\begin{prop}[{Hausdorff-Bernstein-Widder Theorem~\cite[p.~161, Theorem~12b]{widder}}] \label{Widder-thm}
A necessary and sufficient condition that $f(x)$ should be completely monotonic for $0<x<\infty$ is that
\begin{equation} \label{berstein-1}
f(x)=\int_0^\infty e^{-xt}\td\alpha(t),
\end{equation}
where $\alpha(t)$ is non-decreasing and the integral converges for $0<x<\infty$.
\end{prop}

Combining the complete monotonicity in Theorem~\ref{CM-Exp-thm} and the integral representation~\eqref{h(t)-int-rep-eq} with the necessary and sufficient condition in Proposition~\ref{Widder-thm}, it was revealed in~\cite{simp-exp-degree-ext.tex} that
\begin{equation}\label{lca-qi-wang-final-ineq}
\frac{I_1\bigl(2\sqrt{u}\,\bigr)}{\sqrt{u}\,} \ge\frac{u}{1-e^{-u}},\quad u>0.
\end{equation}
Replacing $2\sqrt{u}\,$ by $t$ in~\eqref{lca-qi-wang-final-ineq} yields~\cite[Theorem~1.3]{simp-exp-degree-ext.tex} below.

\begin{theorem}\label{Bessel-2-ineq-thm}
For $t>0$, we have
\begin{equation}\label{I=1-exp-ineq}
I_1(t)>\frac{(t/2)^3}{1-e^{-(t/2)^2}}.
\end{equation}
\end{theorem}

We note that the complete monotonicity in Theorem~\ref{CM-Exp-thm} is the basis of the inequality~\eqref{I=1-exp-ineq} and some results in the subsequent papers~\cite{Bessel-ineq-Dgree-CM.tex, QiBerg.tex}.
\par
The aim of this paper is, with the help of the integral representation~\eqref{h(t)-int-rep-eq} but without using Proposition~\ref{Widder-thm}, to supply a new proof of Theorems~\ref{CM-Exp-thm} and~\ref{Bessel-2-ineq-thm} in a converse direction with that in~\cite{Bessel-ineq-Dgree-CM.tex, QiBerg.tex, simp-exp-degree-ext.tex}. In other words, Theorem~\ref{Bessel-2-ineq-thm} will be firstly and straightforwardly proved, and then Theorem~\ref{CM-Exp-thm} will be done.

\section{A new proof of Theorems~\ref{CM-Exp-thm} and~\ref{Bessel-2-ineq-thm}}

By the definition of the modified Bessel function $I_\nu(z)$ in~\eqref{I=nu(z)-eq}, it is easy to see that
\begin{equation*}
\frac{I_1\bigl(2\sqrt{u}\,\bigr)}{\sqrt{u}\,} =\sum_{k=0}^\infty\frac1{k!\Gamma(k+2)}u^{k} >1+\frac12u+\frac1{12}u^2.
\end{equation*}
Hence, in order to prove~\eqref{lca-qi-wang-final-ineq}, it suffices to show
\begin{equation}\label{Q(u)-dfn}
1+\frac12u+\frac1{12}u^2\ge \frac{u}{1-e^{-u}}
\end{equation}
which is equivalent to
\begin{align*}
&\quad e^u \bigl(12-6 u+u^2\bigr)-12-6 u-u^2\\
&>\biggl(1+u+\frac{u^2}{2} +\frac{u^3}{3!} +\frac{u^4}{4!} +\frac{u^5}{5!}\biggr) \bigl[3+(u-3)^2\bigr]-12-6 u-u^2\\
&=\frac{1}{120} u^5 \biggl[\frac34+\biggl(\frac12-u\biggr)^2\biggr]\\
&\ge0.
\end{align*}
Consequently, the proof of the inequality~\eqref{lca-qi-wang-final-ineq}, that is, Theorem~\ref{Bessel-2-ineq-thm}, is thus complete.
\par
Substituting the inequality~\eqref{lca-qi-wang-final-ineq} into the integral representation~\eqref{h(t)-int-rep-eq} leads to $h(t)>0$ and for $k\in\mathbb{N}$
\begin{equation*}
(-1)^kh^{(k)}(t)=\int_0^\infty \biggl[\frac{I_1\bigl(2\sqrt{u}\,\bigr)}{\sqrt{u}\,} -\frac{u}{1-e^{-u}}\biggr]u^ke^{-tu}\td u>0
\end{equation*}
on $(0,\infty)$. As a result, the function $h(t)$ is completely monotonic on $(0,\infty)$.
\par
The limit~\eqref{h(t)-limit=1} follows immediately from taking $t\to\infty$ on both sides of the integral representation~\eqref{h(t)-int-rep-eq}. Theorem~\ref{CM-Exp-thm} is proved.

\begin{rem}
The inequality~\eqref{Q(u)-dfn} is equivalent to
\begin{equation*}
Q(u)=e^u \bigl(12-6 u+u^2\bigr)-12-6 u-u^2>0.
\end{equation*}
An immediate differentiation yields
\begin{align*}
Q'(u)&=e^u \bigl(u^2-4 u+6\bigr)-2 (u+3),\\
Q''(u)&=e^u \bigl(u^2-2 u+2\bigr)-2,\\
Q'''(u)&=u^2e^u.
\end{align*}
Since $Q'''(u)$ and $Q''(0)=0$, it follows that $Q''(u)>0$ on $(0,\infty)$. Owing to $Q'(0)=0$ and $Q''(u)>0$, it is derived that $Q'(u)>0$. Finally, since $Q(0)=0$, the function $Q(u)$ is positive on $(0,\infty)$. This gives an alternative proof of the inequality~\eqref{Q(u)-dfn}.
\end{rem}

\end{document}